\documentclass[reqno,11pt]{amsart}    
\usepackage{epsf,subfigure}   

\newtheorem{theorem}{Theorem}[section]   
\newtheorem{lemma}[theorem]{Lemma}   
\newtheorem{corollary}[theorem]{Corollary}

\theoremstyle{definition}

\theoremstyle{remark}   
   
\newtheorem{notation}[theorem]{Notation}    
\numberwithin{equation}{section}   
   
\def\graph{{\rm graph\, }}   
\def\R{{\mathbb{R}}}

\def\tru{{\rm tr}\left(u^{ij}\right)}   
\def\trup{\left(1+{\rm tr}\left(u^{ij}\right)\right)}

\def\ol#1{\overline{#1}}  
  
\def\theta{\vartheta}
\def\phi{\varphi}
\def\epsilon{\varepsilon}
\def\e{\varepsilon}
\def\tdot{{\tfrac{\partial}{\partial t}}}
   
\begin{document}    
  
\title{Translating solutions for Gau{\ss} curvature flows
  with Neumann boundary conditions}  
   
\author{Oliver C.\ Schn\"urer}   
\address{Max Planck Institute for Mathematics in the Sciences,    
  Inselstr.\ 22-26, 04103 Leipzig, Germany}
\email{Oliver.Schnuerer@mis.mpg.de}   
   
\author{Hartmut R.\ Schwetlick}   
\address{Max Planck Institute for Mathematics in the Sciences,   
  Inselstr.\ 22-26, 04103 Leipzig, Germany}   
\email{Hartmut.Schwetlick@mis.mpg.de}   
%\thanks{}   
   
%    General info   
\renewcommand{\subjclassname}{%  
  \textup{2000} Mathematics Subject Classification}  
   
\date{August 2002, Revised February 2003.}   
   
%\dedicatory{}   
   
%\keywords{}   
   
\begin{abstract}   
We consider strictly convex hypersurfaces
which are evolving  
by the non-parametric logarithmic Gau{\ss} 
curvature flow subject to a Neumann boundary condition.
Solutions are shown to converge smoothly to hypersurfaces moving
by translation. In particular, for bounded 
domains we prove that convex functions with prescribed normal
derivative satisfy a uniform oscillation estimate.
\end{abstract}   
\maketitle   
   
\markboth{Oliver C.\ Schn\"urer and Hartmut R.\ Schwetlick}  
{Translating solutions for Gau{\ss} curvature flows
with Neumann boundary conditions}  

\section{Introduction}

In this paper, we evolve hypersurfaces represented as
graphs of strictly convex functions over  strictly convex bounded domains
by the non-parametric logarithmic Gau{\ss} curvature
flow subject to a Neumann boundary condition. We show that 
solutions exist for all time and converge smoothly to 
translating solutions.

To be more precise, we address the following slightly
more general problem: Let $u_0$ be a strictly convex
function over a smooth strictly convex bounded domain
$\Omega\subset\R^n$. 
We use the phrase {\it strictly convex\/} for functions
whose Hessian is positive definite, and for domains for which
all principal curvatures of the boundary are positive. 
Assume that $u_0$ is smooth up to the boundary,
$u_0\in C^\infty\left(\ol\Omega\right)$, and satisfies
$$D_\nu u_0=\phi\quad\text{on~}\partial\Omega,$$
where $\nu$ is the inner unit normal to $\partial\Omega$
and $\phi\in C^\infty(\partial\Omega)$. 
Let $f\in C^\infty\left(\ol\Omega\times\R^n\right)$.
We prove the following
\begin{theorem}\label{main thm}
For $\Omega$, $\phi$, $\nu$, $f$
and $u_0$ as introduced above,
there exists a family $u(\cdot,t)$, $t\in[0,\infty)$, of strictly convex
functions solving
\begin{equation}\label{flow eqn}
\left\{\begin{array}{r@{~=~}ll}
\tdot u & \log\det D^2u-\log f(x,Du) & \mbox{in~}\Omega\times
[0,\infty),\\
D_\nu u (\cdot, t)& \phi & 
\mbox{on~}\partial\Omega,\,t>0,\\
u(\cdot,0) & u_0 & \mbox{in~}\Omega,
\end{array}\right.
\end{equation}
where $u\in C^\infty\left(\ol\Omega\times(0,\infty)\right)$,
and  $u(\cdot,t)$ approaches $u_0$ in $C^2\left(\ol\Omega\right)$ as $t\to0$.
Moreover, $u(\cdot,t)$ converges
smoothly to a translating solution, i.e., to
a solution with constant time derivative.  
\end{theorem}
We remark that the parabolic maximum principle implies that
the asymptotic solutions for different initial data $u_0$ 
are unique up to a constant. 
For the Gau{\ss} curvature flow, mentioned in the beginning,
the flow equation takes the form
$$\tdot u=\log\frac{\det D^2u}{\left(1+|Du|^2\right)^{\frac{n+2}{2}}}=
\log\det D^2u-\tfrac{n+2}{2}\log\left(1+|Du|^2\right).$$

In the proof of Theorem \ref{main thm}, we generalize 
the result to
the so-called oblique boundary condition
$$D_\beta u=\phi\quad\text{on~}\partial\Omega,$$
where $\beta$ is a unit vector field which is $C^1$-close to $\nu$, 
i.e., such that there exists
a small positive constant $c_\beta>0$ for which
$\Vert\nu-\beta\Vert_{C^1}\le c_\beta$.
Such a generalization
of a Neumann boundary condition is studied for the elliptic case
in \cite{Urbas Calc 1998}.

We base the barrier construction in the proof of Theorem \ref{main thm}
on solutions to a related elliptic problem given by the following
\begin{theorem}\label{ell thm}
Consider $\Omega$, $\nu$, $\phi$, and $f$ as introduced before.
There exist a unique $v\in\R$ and a strictly convex
function $u\in C^\infty\left(\ol\Omega\right)$ solving
the boundary value problem 
\begin{equation}\label{ell cor eqn}
\left\{\begin{array}{r@{~=~}ll}
\det D^2u&e^v\cdot f(x,Du) & \mbox{in~}\Omega,\\
D_\nu u & \phi & \mbox{on~}\partial\Omega,
\end{array}\right.
\end{equation}
provided that 
there exists a smooth strictly convex function $u_0$ 
satisfying the boundary condition 
$D_\nu u_0=\phi$. The function $u$ is unique up to a 
constant.
\end{theorem}
We remark that translating solutions to (\ref{flow eqn}) can be viewed
as solutions to  (\ref{ell cor eqn}), where $v$ denotes the speed.

A situation similar
to Theorem \ref{main thm} is considered for the mean 
curvature flow in \cite{hui bdy}. Hypersurfaces of prescribed
Gau{\ss} curvature subject to Neumann boundary conditions
are found in the pioneering paper \cite{ltu}.
The extension to the oblique boundary value problem is made 
in \cite{Urbas Calc 1998}.
Flows of Monge-Amp\`ere type for the Neumann and the second boundary
value problem are studied in \cite{osks2001}.
In our setting, 
the situation is more degenerate as neither $f$ nor $\phi$ do depend
on $u$. Thus, both $f$ and $\phi$ fail to satisfy
the crucial monotonicity requirement with respect to $u$.
For the second boundary value problem, translating
solutions to flows of Gau{\ss} curvature type
are considered in \cite{os2001trans}.

As mentioned in \cite{os2001trans},
methods of \cite{osks2001} can be adapted 
to  the non-parametric
logarithmic Gau{\ss} curvature flow
subject to the second boundary condition. 
For the Neumann boundary value
problem, however, the lack of monotonicity requires
a new proof to  uniformly bound  the oscillation of a solution.
Therefore, we establish a generalization
of the spatial $C^1$-estimates of \cite{ltu}. 
Then, we use the translating solutions provided by Theorem \ref{ell thm},
in particular its uniquely determined speed, to construct an
auxiliary barrier function and to obtain
uniform spatial $C^2$-estimates.
Hence, the results
of Krylov, Safonov, Evans and Schauder
imply uniform bounds on higher derivatives
for all times, uniformly bounded away from $0$. 
This allows to show smooth convergence to a translating solution.

The paper is organized as follows. 
As explained in \cite{osks2001}, standard linear parabolic theory 
and the implicit function theorem imply
short-time existence. 
We show uniform first-order estimates in Sections 2 and 3.
Section 4 contains the proof of  Theorem \ref{ell thm}.
Having the unique velocity of a translating solution, we can 
prove uniform a priori $C^2$-estimates in Section 5. 
Finally, in Section 6, we prove smooth convergence to a translating solution.
In Appendix A we apply  Theorem \ref{ell thm} to construct entire graphs of
prescribed Gau{\ss} curvature.
To illustrate the convergence of the flow, we carry out a numerical
integration on a planar domain in Appendix B.

We wish to thank J\"urgen Jost, Stefan M\"uller and the 
Max Planck Institute for Mathematics
in the Sciences for their hospitality.
The first author also wishes to acknowledge that this paper was
  finished while he was at Harvard University, supported by the
  Humboldt foundation. 

\section{$\tdot u$-Estimate}
\begin{notation}
We write a dot to denote the time derivative, whereas
we use indices for the spatial partial derivatives.
Let $f_{p_i}$ denote a derivative of $f$ with respect
to the gradient. 
For a vector  $\xi$ we define $u_\xi\equiv\xi^iu_i$.
For the logarithm of $f$ we
use $\hat f\equiv\log f$. We use the Einstein
summation convention and sum over repeated
upper and lower indices. The inverse of the
Hessian of $u$ is denoted by $\left(u^{ij}\right)=(u_{ij})^{-1}$.
We remark that -- besides in the case $u^{ij}$ -- indices are
lifted with respect to the Euclidean metric.
The letter $c$ denotes a generic positive constant.
Furthermore, we may assume that
$0\in\Omega$.
\end{notation}

\begin{lemma}\label{u dot lem}
Under the assumptions of Theorem \ref{main thm}, there holds
$$|\dot u|\le\max_{t=0}|\dot u|,$$
as long as a smooth convex solution of (\ref{flow eqn}) 
exists.
\end{lemma}
\begin{proof}
Similar to \cite{osks2001}, we consider 
$$r:=(\dot u)^2.$$
We get the evolution equation
$$\dot r=u^{ij}r_{ij}-2u^{ij}\dot u_i\dot u_j  
-\hat f_{p_i}r_i.$$
Hence, $\dot r\le0$ at a maximum 
of $r$ in $\Omega\times[0,t]$. 
Now, 
assume that a maximum of $r$ on
$\ol\Omega\times[0,t]$ occurs 
at $(x_0,t_0)$ with 
$x_0\in\partial\Omega$. 
If $t_0=0$ the lemma holds. 
Thus, in the following, we may assume $t_0>0$.
If $r$ is constant,
then $u$ is a translating solution, and our
lemma holds. Otherwise, we  get $r_\beta(x_0)<0$ from
the Hopf boundary point lemma.
At $x_0$ we compute
$$0>r_\beta=\left(\left(\dot u\right)^2\right)_\beta
=2\dot u\dot u_\beta=2\dot u\tdot\phi(x)=0.$$
Contradiction.
Note, that the assumption  $t_0>0$ guarantees that $u$ is smooth near
 $(x_0,t_0)$ allowing to 
interchange differentiation with respect to time and space.
\end{proof}

Integrating the last estimate yields
\begin{corollary}
As long as a smooth convex solution of (\ref{flow eqn})
exists, we obtain the estimate
$$|u(x,t)|\le\sup\limits_\Omega|u_0|
+\sup\limits_\Omega|\dot u(\cdot, 0)|\cdot t.$$
\end{corollary}

\section{Ice-cream cone estimate}
 The following theorem generalizes the 
$C^1$-estimates of \cite{ltu}.
It is essential for our situation, because the oscillation but not the
$C^0$-norm of the solution is expected to be bounded uniformly in time.

We wish to mention oscillation estimates of Urbas 
\cite{Urbas JDG 1984, Urbas xxx}.
There, the convexity of $u$ and appropriate growth of $f(x,p)$ in $p$
is used, whereas Theorem 
\ref{iccl} combines convexity and the boundary condition.

\begin{theorem}[Ice-cream cone estimate]\label{iccl}
Let $\Omega\subset\R^n$ be a smooth bounded domain,
$u:\ol\Omega\to\R$ a smooth strictly convex function
with $|u_\beta|$ uniformly bounded on $\partial\Omega$,
where $\beta$ is a unit vector field on $\partial\Omega$ such
that  $\langle\beta,\nu\rangle\ge \tilde c_\beta$
for a positive constant $\tilde c_\beta>0$
(recall that $\nu$ is the inner unit normal to $\partial\Omega$). 
Then there is a uniform bound for $\sup|Du|$, independent of $\sup|u|$.
\end{theorem}
In view of  Lemma \ref{u dot lem}, the result above yields an estimate for 
the
full $C^1$-norm of solutions $u$ to (\ref{flow eqn}).
Note, that only the estimate for the derivatives of $u$ 
is uniform in time.

\begin{proof}
In the name of the theorem we want to emphasize that our proof uses balls
and cones, similar to ice-cream placed in a cone of waffle.
We argue by contradiction.
Assume that
there exists a point $x_0$, where
$|Du|$ is maximal and equal to $M$.
If $M$ is larger than a suitably chosen constant $M_0$
we will find a contradiction.
As $u$ is strictly convex, we see that $x_0\in\partial\Omega$.
At $x_0$ we find a tangential direction $\xi_0$ such that
$\langle Du(x_0),\xi_0\rangle$ is maximal compared to 
all other tangential directions. 
Here and later, unit vectors are called directions.
We wish to prove a lower estimate for $\langle Du(x_0),\xi_0\rangle$
in terms of $M$.

Let $\xi_1$ be a direction such that $\langle Du(x_0),\xi_1\rangle
=M$. 
Similar to \cite{Urbas Calc 1998}, we decompose a direction $\xi$ using $\beta$ and a
tangential vector $\tau(\xi)$ as
\begin{equation}\label{dec xi}
\xi=\tau(\xi)+\frac{\langle\nu,\xi\rangle}{\langle\beta,\nu
\rangle}\beta,
\end{equation}
where
$$\tau(\xi)=\xi-\langle\nu,\xi\rangle\nu
-\frac{\langle\nu,\xi\rangle}{\langle\beta,\nu\rangle}\beta^T,
\quad\beta^T=\beta-\langle\beta,\nu\rangle\nu.$$
Note, that $|\tau(\xi)|$ is bounded by assumption. 
Decomposing  $\xi_1$, we get
\begin{eqnarray*}
M&=&\langle Du,\xi_1\rangle=\langle Du,\tau(\xi_1)\rangle
+\frac{\langle\nu,\xi_1\rangle}{\langle\beta,\nu\rangle}
\langle Du,\beta\rangle\\
&\le&|\tau(\xi_1)|\cdot\max\limits_{\genfrac{}{}{0pt}{}
{\tau\in T_{x_0}\partial\Omega}{|\tau|=1}}
\langle Du,\tau\rangle+c\\
&=&|\tau(\xi_1)|\cdot\langle Du(x_0),\xi_0\rangle+c.
\end{eqnarray*}
Hence, we
deduce that $\langle Du(x_0),\xi_0\rangle\ge\frac{M}{c}$,
as long as $M\ge M_0$ and 
$M_0$ is chosen sufficiently large. 
For a direction $\xi$ near $\xi_0$,
say  $|\xi-\xi_0|<\e=\frac1{2c}<1$, we obtain
\begin{equation}\label{cone found}
\langle Du(x_0),\xi\rangle=
\langle Du,\xi_0\rangle+\langle Du,\xi-\xi_0\rangle
\ge\genfrac{}{}{}1{M}{c}-M|\xi-\xi_0|\ge\e M.
\end{equation} 
 From the convexity of $u$, we deduce that
$\langle Du(y),\xi\rangle\ge\e M$
for all points $y\in\ol\Omega$ of the form 
$y=x_0+\lambda\cdot\xi$.
Here, $\lambda>0$ and $\xi$ are chosen such that  
$|\xi-\xi_0|\le\epsilon$ and
$x_0+t\cdot\lambda\cdot\xi\in\ol\Omega$ for all $t\in[0,1]$.

The uniform boundedness
of the principal curvatures of $\partial\Omega\subset\R^{n}$ implies
that there exist $R>0$ and $x_1\in\partial\Omega$
such that  $|x_0-x_1|>2R$, and, especially,
any $x\in B_R(x_1)\cap\partial\Omega$ can be written in the form 
$x_0+\lambda\cdot\xi$,
as described above.
Thus, according to  (\ref{cone found}),
 $|Du|\ge\e M$ in
$\partial\Omega\cap B_{R}(x_1)$.
Due to our construction, we have
$$\inf\limits_{x\in B_R(x_1)\cap\partial\Omega}u(x)>u(x_0).$$
\begin{figure}[htb]
\typeout{includes ice2.eps}
\epsfxsize=0.9\textwidth
\centerline{\epsfbox{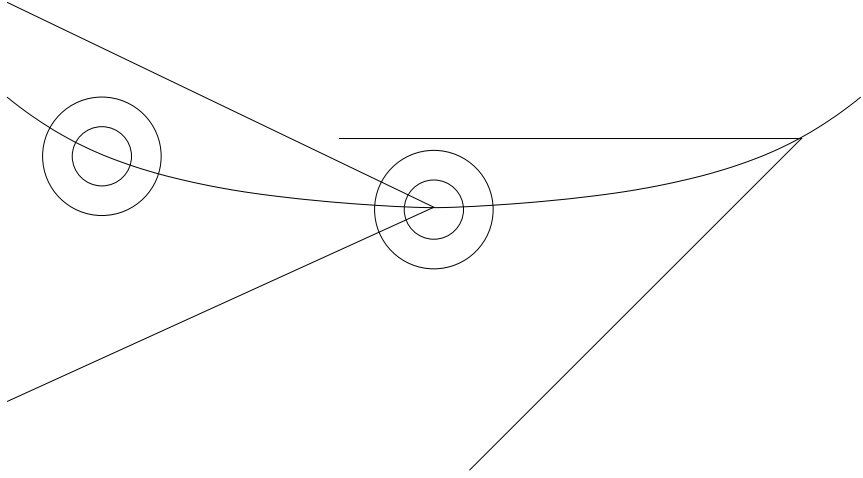}}
\caption{Ice-cream cone estimate} 
\label{ice pic}
\end{figure}%
{}Figure \ref{ice pic} shows a part of $\partial\Omega$
and two cones corresponding to the directions $\xi$
as well as two pairs of concentric balls.
The larger ones are the balls $B_R$ mentioned above,
the smaller ones are  introduced in the following.

Now, we  proceed iteratively.
Note, that $R$ and $\varepsilon$ can be chosen as fixed
constants, independent of the point $x_0$.
As long as  $|Du(x_i)|\ge M \e^i\ge M_0$,
we can find a further point $x_{i+1}$ 
as we went from $x_0$ to $x_1$.
Thus, for $M=\sup|Du|$  sufficiently large,
 we can construct 
a sequence of points $\{x_i\}_{i=0,\ldots,N}$ of arbitrarily large length $N$,
satisfying for all $i\ge1$
$$
|Du|\ge M\e^i\quad\text{on~}\partial\Omega\cap B_{R}(x_i)
$$
and
$$
\inf\limits_{x\in B_R(x_i)\cap\partial\Omega}u(x)>u(x_{i-1}).
$$
Since $\partial\Omega$ has finite measure and bounded principal curvatures,
 there is an
upper bound  $N_0(\rho)$ on the number of pairwise disjoint restricted balls
$B_\rho(y_j)\cap\partial\Omega$ for fixed $\rho>0$ and $y_j\in\partial\Omega$.

Hence, if $M=\sup|Du|>M_0\e^{-N_0\left(\frac R2\right)}$, 
there will be two points $x_{i_0},x_{j_0}$ with $i_0>j_0>0$
such that 
$$B_{\frac R2}(x_{i_0})\cap B_{\frac R2}(x_{j_0})
\cap\partial\Omega\not=\emptyset.$$
But $x_{j_0}\in B_R(x_{i_0})$ implies
$$u(x_{j_0})<u(x_{j_0+1})<\ldots<u(x_{i_0-1})<
\inf\limits_{x\in B_R(x_{i_0})\cap\partial\Omega}u(x)
\le u(x_{j_0}).$$ 
Contradiction.
\end{proof}

\section{Existence of translating solutions}
\label{ell sec}

This section is devoted to the proof of Theorem \ref{ell thm}.
That is, we construct
solutions to the elliptic problem
\begin{align}\label{ell eqn}
\left\{\begin{array}{r@{~=~}ll}
\displaystyle 
v&\log\det D^2u-\log f(x,Du)\quad&\text{in }\Omega,\\
\displaystyle 
u_\beta&\phi&\text{on }\partial\Omega.
\end{array}\right.
\end{align}
The absence of any monotonicity property in $u$, 
both in  $f$ as well as
in the boundary condition $\phi$, 
limits seriously the existence of solutions.

Step 1:
Here, we show that for given $\e>0$ and $v\in\R$
there is a  unique solution $u_{\e,v}$ of
\begin{align*}
\left\{\begin{array}{r@{~=~}ll}
\displaystyle 
\det D^2u&e^vf(x,Du)e^{\e u}\quad&\text{in }\Omega,\\
\displaystyle 
u_\beta&\phi&\text{on }\partial\Omega.
\end{array}\right.
\tag{$*_{\e,v}$}
\end{align*}
Note, that the dependence on $v$ is continuous and strictly decreasing.
In fact, we have the explicit relation
\[
u_{\e,v}=u_{\e,0}-\tfrac v\e.
\]
To show the unique existence of $u_{\e,v}$, we will derive 
an  a priori $C^0$-bound, 
then the ice-cream cone estimate, Theorem \ref{iccl}, yields the $C^1$-bound.
Having  controlled the full $C^1$-norm, we can estimate the $C^2$-norm
exactly as in  Urbas \cite{Urbas Calc 1998}, since the strict  monotonicity
assumption on $\phi$ is not used for this part.
A detailed argument, applied to  the parabolic case,
is given in the next section.
Bounds for higher $C^k$-norms follow via the estimates due to 
Krylov, Safonov, Evans, and from Schauder theory.

To get the $C^0$-bound, we define suitable barriers for ($*_{\e,0}$).
Recall, that $u_0$ is a convex function satisfying 
$(u_0)_\beta=\phi$ on $\partial\Omega$.
Define $u_\e^\pm=u_0\pm M/\e$, where $M>0$ will be chosen later.
Then 
\[
\frac{\det D^2u_\e^\pm}{f(x,Du_\e^\pm)e^{\e u_\e^\pm}}=
\frac{\det D^2u_0}{f(x,Du_0)e^{\e u_0\pm M}}=
g(x)e^{-\e u_0}e^{\mp M},
\]
with $c^{-1}<g(x)<c$.
Hence, restricting ourselves to $\e<1$, 
there exists a large constant $M>0$, not depending on $\e$,
such that  
$u_\e^+$ is a strict supersolution and
$u_\e^-$ is a strict subsolution of $(*_{\e,0})$.
This implies, that
\begin{equation}\label{ue0 eqn}
u_\e^-< u_{\e,0}<u_\e^+ ,
\end{equation}
or equivalently,
\[
\left|u_{\e,v}-(u_0-\tfrac v\e)\right|<\tfrac M\e.
\]

Step 2: Now, we consider the limit $\e\to0$.
In general, we cannot expect that the sup bounds for $u_{\e,v}$ can be
obtained uniformly in $\e$. 
In fact, it follows from the maximum principle  that only for 
a unique $v$, there is a solution to
($*_{\e,v}$) with $\e=0$. 
Observe, that (\ref{ue0 eqn}) implies that
\[
u_{\e,+M}< u_0< u_{\e,-M}.
\]
Therefore,  for every $\e>0$ we can find a unique $v_\e\in(-M,M)$ such that
$u_{\e,v_\e}(0)=u_0(0)$. 
Note, that  (\ref{ue0 eqn}) does not suffice to control the oscillation of
$u_{\e,v_\e}$,  uniformly in $\e$. We employ
the ice-cream cone estimate
to bound $u_{\e,v_\e}$  uniformly  in $C^1$.
Again, uniform $C^1$-bounds imply uniform higher $C^k$-bounds.

Now, we choose a sequence $\e_i\to0$ as $i\to\infty$.
Since $v_{\e_i}$ is bounded, there exists a subsequence, relabeled,
such that $v_{\e_i}\to v^\infty$ and 
$u_{\e_i,v_{\e_i}}\to u^\infty_{\text{ell}}$ in any  $C^k$-norm.
This completes the proof of Theorem \ref{ell thm}.

The extension $u^\infty(x,t):=u^\infty_{\text{ell}}(x)+v^\infty t$ 
is a translating solution, as,  by construction,
  $u^\infty$ satisfies
$$\left\{\begin{array}{r@{~=~}ll}
\displaystyle v^\infty=\dot u^\infty&\log\det D^2u^\infty
-\log f(x,Du^\infty)&\mbox{in }\Omega,\\
\displaystyle u^\infty_\beta&\phi&\mbox{on }\partial\Omega.
\end{array}\right.$$
     
\section{Parabolic $C^2$-estimates}

The following argument is a modification of the proofs
in \cite{ltu, osks2001, Urbas Calc 1998}.
We use the translating solution $u^\infty$, especially,
its speed $v^\infty$, to construct an auxiliary barrier
function. 

Assume that $u^\infty_{\text{ell}}>u_0$.  We define 
$$\tilde\phi(x,z)=\phi(x)+(z-u^\infty_{\text{ell}}).$$
Due to uniform estimates on the gradient
of $u$, we can find a positive constant $\mu_0$ such that 
$$\min\{f(x,Du),\,f(x,Du^\infty_{\text{ell}})\}\ge\mu_0.$$
For $0<\rho<1$, consider the elliptic boundary value problem 
\begin{equation}\label{int psi}
\left\{\begin{array}{r@{~=~}ll}
v^\infty & \log\det D^2\psi-\log\frac{\mu_0}{2} & 
\mbox{in~}\Omega,\\
\psi_\beta & \tilde\phi\left(x,\psi+\rho\cdot\vert x\vert^2
\right)-2\rho\,\langle x,\beta\rangle & \mbox{on~}
\partial\Omega.
\end{array}\right.
\end{equation}
We wish to show a uniform a priori $C^2$-estimate for $\psi$.
Theorem \ref{iccl} gives an estimate for the gradient.
Similarly to \cite{Urbas Calc 1998}, bounds on the second
derivatives follow.
It is here that the smallness of $\Vert\nu-\beta\Vert_{C^1}$
is used.
Thus, it remains to prove uniform 
$C^0$-estimates. 
Note, that a convex solution $\psi$ cannot satisfy
$\psi_\beta(x)>0$ for all $x\in\partial\Omega$.
Hence, the upper bound on $\psi$ follows since  $\tilde\phi(x,z)
\to\infty$ uniformly as $z\to\infty$. For the lower bound,
we consider $\psi-u^\infty_{\text{ell}}$. Applying the maximum principle
to the differential inequality
$$\left\{\begin{array}{rcll}
0&>&\log\det D^2\psi-\log\det D^2u^\infty_{\text{ell}}&\mbox{in~}\Omega,\\
(\psi-u^\infty_{\text{ell}})_\beta&
=&\psi-u^\infty_{\text{ell}}+\rho\cdot|x|^2-2\rho\,\langle
x,\beta\rangle &\mbox{on~}\partial\Omega,
\end{array}\right.$$
we see that
$\psi-u^\infty$ cannot attain an interior minimum. If
a minimum occurs on $\partial\Omega$, we get
$$0\le(\psi-u^\infty_{\text{ell}})_\beta
=\psi-u^\infty_{\text{ell}}+\rho\cdot|x|^2
-2\rho\,\langle x,\beta\rangle.$$
Thus, $\psi$ is uniformly bounded below,
a solution to (\ref{int psi}) exists. 
Due to the uniform $C^2$-estimates, we can
fix $\lambda>0$ such that 
\begin{equation}\label{lambda fix}
 (\psi_{ij})\ge\lambda\rm{Id}.
\end{equation}
Furthermore, these estimates allow
to fix $\rho>0$ such that $\ol\psi:=\psi+\rho\cdot\vert x\vert^2$
satisfies
$$\left\{\begin{array}{rcll}
v^\infty&>&\log\det D^2\ol\psi-\log\mu_0&\mbox{in~}\Omega,\\
\ol\psi_\beta&=&\tilde\phi\left(x,\ol\psi\right)&
\mbox{on~}\partial\Omega.
\end{array}\right.$$
Applying the maximum principle, we get 
$u^\infty_{\text{ell}}\le\ol\psi$.
We extend $\psi$ and $\ol\psi$ by setting
$\psi(x,t):=\psi(x)+t\cdot v^\infty$,
$\ol\psi(x,t):=\ol\psi(x)+t\cdot v^\infty$,
respectively. Thus, 
$u\le u^\infty\le\ol\psi$,
where $u^\infty$ is the translating solution defined 
in Section \ref{ell sec}. We get for  $x\in\partial\Omega$
$$\left(\ol\psi_\beta-u_\beta\right)(x,t)
=\ol\psi_\beta(x,0)-\phi(x)
=\left(\ol\psi-u^\infty\right)(x,0)=
(\ol\psi-u^\infty)(x)\ge0.$$
Furthermore, for a sufficiently small $\delta_0>0$,
\begin{equation}\label{int delta0}
(\psi-u)_\beta=\left(\ol\psi-\rho\cdot|x|^2-u\right)_\beta
\ge-2\rho\,\langle x,\beta\rangle\ge\delta_0>0
\quad\text{on~}\partial\Omega,
\end{equation}
provided that $\beta$ is $C^0$-close to $\nu$. Here, we used
that $0\in\Omega$ implying $\langle x,\nu\rangle<0$.
Using these preparations, we prove a priori $C^2$-%
estimates similarly to \cite{osks2001, Urbas
Calc 1998}. For the
reader's convenience, we repeat the arguments
incorporating the necessary modifications to the parabolic case.

\subsection{Preliminary results}
Assume, that a smooth solution $u$ of our flow equation (\ref{flow eqn})
exists on the time interval $[0,T]$. 
We will use the letter $\tau$ to indicate a direction tangential 
to $\partial\Omega$.

\begin{lemma}[Mixed $C^2$-estimates at the boundary]  
\label{mixed lemma}  
Let $u$ be a solution of (\ref{flow eqn}).
Then $|u_{\tau\beta}|$ remains uniformly bounded
on $\partial\Omega$.
\end{lemma}
\begin{proof}
We represent $\partial\Omega$ locally as $\graph\omega$ over  
its tangent plane at a fixed point $x_0\in\partial\Omega$
such that locally  
$\Omega=\{(\hat x, x^n):x^n>\omega(\hat x)\}$.
Let us extend $\beta$ and $\phi$ smoothly.
At $x_0$, 
differentiating the oblique boundary condition
$$\beta^i(\hat x)u_i(\hat x,\omega(\hat x))=
\varphi(\hat x,\omega(\hat x)),
\quad\hat x\in\R^{n-1},$$
with respect to tangential directions 
$\hat x^j$, $1\le j\le n-1$,
$$\beta^i_ju_i+\beta^iu_{ij}+\beta^iu_{in}\omega_j=
\varphi_j+\varphi_n\omega_j,$$
we obtain at $x_0\equiv(\hat x_0,\omega(\hat x_0))\in\partial  
\Omega$ a bound for
$\beta^iu_{ij}$. Thereby, we use the gradient estimate
for $u$ and
$D\omega(\hat x_0)=0$.
Multiplying with $\tau^j$ gives the result.  
\end{proof}

\begin{lemma}[Double oblique $C^2$-estimates at the boundary]  
\label{Double normal lemma}
 For any solution of (\ref{flow eqn}),
$|u_{\beta\beta}|$ is uniformly bounded on $\partial\Omega$.
\end{lemma}
\begin{proof}
Note that $u_{\beta\beta}>0$ as $u(\cdot,t)$ is strictly 
convex for each $t$. 
We keep the geometric  
setting of the proof of Lemma \ref{mixed lemma}  
with $x_0\in\partial\Omega$.
 From (\ref{flow eqn}) we obtain
$$\dot u_k=u^{ij}u_{ijk}-(\hat f_k+\hat f_{p_i}u_{ik}).$$  
We define 
$$Lw:=\dot w-u^{ij}w_{ij}+\hat f_{p_i}w_i.$$
We can find   
appropriate extensions of $\beta$ and $\varphi$
such that
$$\left|L\left(\beta^ku_k-\varphi(x)\right)\right|\le c\cdot\trup.$$
We choose $\delta>0$ sufficiently small and define
$\Omega_\delta:=\Omega\cap B_\delta(x_0)$. Set
$$\vartheta:=d-\mu d^2,$$
where $\mu\gg1$ is chosen sufficiently large, and 
$d$ denotes the distance
to $\partial\Omega$. We will show that
in $\Omega_\delta$ there holds $L\vartheta\ge
\frac{\varepsilon}{3}\tru$ for a small constant $\varepsilon>0$,
depending only on a positive lower bound for the principal  
curvatures of $\partial\Omega$. Next,
\begin{eqnarray}\nonumber
L\vartheta&=&-u^{ij}d_{ij}+2\mu u^{ij}d_id_j
  +2\mu u^{ij}dd_{ij}+\hat f_{p_i}(d_i-2\mu dd_i)\nonumber\\
&\ge&-u^{ij}d_{ij}+2\mu u^{ij}d_id_j
  -c\mu d\trup-c.\nonumber
\end{eqnarray}
We use the strict convexity of $\partial\Omega$, $|Dd-e_n|\le c\delta$,
$|u^{kl}|\le\tru$, $1\le k,\, l\le n$, and the
inequality for arithmetic and geometric means
\begin{eqnarray*}
L\vartheta&\ge&\varepsilon\,\tru+\mu u^{nn}
  -c\mu\delta\trup-c\\
&\ge&\tfrac{n}{3}\left(\det\left(u^{ij}\right)\right)^\frac{1}{n}
  \cdot\varepsilon^\frac{n-1}{n}\cdot\mu^\frac{1}{n}
  +\tfrac{2}{3}\varepsilon\,\tru\\  
&&{}-c\mu\delta\trup-c.
\end{eqnarray*}
More precisely, we used
$$\frac{\epsilon\,\tru+\mu u^{nn}}{3}\ge
\frac{n}{3}\epsilon^{\frac{n-1}{n}}\mu^{\frac{1}{n}}
\left(\prod_{i=1}^{n}u^{ii}\right)^{\frac1n},$$
and, assuming $(u^{ij})_{i,j<n}$ is diagonal,
\begin{eqnarray}
\det\left(u^{ij}\right) &=&
\det\left(\begin{array}{ccccc}
u^{11} & 0       & \cdots & 0             & u^{1n}   \\
0      & \ddots  & \ddots & \vdots        & \vdots   \\
\vdots & \ddots  & \ddots & 0             & \vdots   \\
0      & \cdots  & 0      & u^{n-1\, n-1} & u^{n-1\, n}\\
u^{1n} & \cdots  & \cdots & u^{n-1\, n}   & u^{nn}   \\
\end{array}\right)\nonumber\\[.3em]
&=& \prod_{i=1}^n u^{ii} \:-\: \sum_{i<n} |u^{ni}|^2 
\:
\prod_{\genfrac{}{}{0pt}{}{j\neq i}{j<n}} u^{jj} 
\;\leq\; 
\prod_{i=1}^n u^{ii} \;.
\label{det exp}
\end{eqnarray}

Since 
$$\det\left(u^{ij}\right)=\frac1{\det (u_{ij})}
  =\exp\left(-\hat f-\dot u\right),$$
$\det\left(u^{ij}\right)$ is uniformly bounded from below by  
a positive constant.
  We may choose $\mu$ so large, that 
$$\tfrac n3\left(\det\left(u^{ij}\right)\right)^{\frac1n}
\cdot\epsilon^{\frac{n-1}n}\cdot\mu^{\frac1n}\ge c+1.$$
 For $\delta\le
\frac{1}{c\mu}\min\left\{1,\frac{1}{3}\varepsilon\right\}$,
we get
$$L\vartheta\ge\tfrac{1}{3}\varepsilon\,\tru.$$  
 Furthermore, $\vartheta\ge 0$ on $\partial\Omega_\delta$,
if $\delta$ is chosen smaller if necessary.
\par
Let $l$ be an affine linear function such that $l(x_0)=0$ and 
$$l\ge\beta^i(u_0)_i-\phi
\quad\text{in~}\Omega_\delta.$$

 For constants $A$, $B>0$, consider the function
$$\Theta:=A\vartheta+B|x-x_0|^2-\left(\beta^iu_i-\varphi(x)
\right)+l.$$
We fix $B\gg 1$, get $\Theta\ge0$ on $\partial\Omega_\delta$,
and deduce for $A\gg B$ that $L\Theta\ge0$, since $\tru$ is
bounded from below by a positive constant.
The maximum principle yields $\Theta\ge0$
in $\Omega_\delta$. As $\Theta(x_0)=0$, we 
conclude that $\Theta_\beta(x_0)\ge0$
implying $u_{\beta\beta}\le c$.  
\end{proof}
\begin{lemma}\label{intro lem}
For a solution of (\ref{flow eqn}), there holds
$$\min\limits_{\genfrac{}{}{0pt}{}{t\in[0,T]}{x\in\partial\Omega}}\;
\max\limits_{\genfrac{}{}{0pt}{}{\xi\in T_x\partial\Omega}{|\xi|=1}}
u_{\xi\xi}(x,t)>0.$$
\end{lemma}
\begin{proof}
We have already seen that there is a uniform positive
lower bound for $\det D^2u$. Again, at a fixed boundary point,
we may choose a coordinate system, such that $e_n$ is
equal to the inner unit normal of $\partial\Omega$ 
and $(u_{ij})_{i,\,j<n}$ is diagonal. 
Similarly to (\ref{det exp}), we estimate
\begin{equation}
\det\left(u_{ij}\right)
=\prod_{i=1}^n u_{ii} \:-\: \sum_{i<n} |u_{ni}|^2 
\:\prod_{
\genfrac{}{}{0pt}{}{j\neq i}{j<n}
} u_{jj} \;\leq\; \prod_{i=1}^n u_{ii} \;.
\label{det exp1}
\end{equation}
We decompose $\nu$ as
$$\nu=\frac{1}{\langle\beta,\nu\rangle}
\left(\beta-\beta^T\right),$$
$\beta^T$ as in the proof of Theorem \ref{iccl},
and get in view of 
Lemmata \ref{mixed lemma} and \ref{Double normal
lemma}
$$u_{\nu\nu}\le\frac{1}{\langle\beta,\nu\rangle^2}
\left(\left\vert\beta^T\right\vert^2\cdot
\max\limits_{i<n}u_{ii}+c\right).$$
Finally, the claimed bound follows from (\ref{det exp1}).
\end{proof}

\subsection{Remaining $C^2$-estimates}

Similarly to \cite{Urbas Calc 1998},
we define 
$$w(x,\xi,t):=e^{\alpha(\psi-u)+\gamma\cdot|Du|^2}
\cdot u_{\xi\xi}$$
for $(x,\xi,t)\in\ol\Omega\times S^{n-1}\times[0,T]$
and positive constants $\alpha$, $\gamma$ to be fixed later. 
We assume that $w$, restricted to boundary points 
and tangential directions, attains its 
maximum at $x_w\in\partial\Omega$
in a tangential direction which we may take to be $e_1$,
and $t_w\in[0,T]$. We may assume that $t_w>0$.
Furthermore,
we fix Euclidean coordinates such that 
$e_n$ is the inner normal direction and 
$(u_{ij})_{i,\,j<n}(x_w)$ is diagonal. 
Decompose $e_1$ as
$$e_1=\tau(e_1)+\frac{\langle\nu,e_1\rangle}
{\langle\beta,\nu\rangle}\beta,$$
where
$$\tau(e_1)=\tau=e_1-\langle\nu,e_1\rangle\nu
-\frac{\langle\nu,e_1\rangle}{\langle\beta,
\nu\rangle}\beta^T,\quad\beta^T=\beta-\langle\beta,\nu\rangle
\nu.$$
Note that $\tau$ is tangential, but not necessarily of 
unit length. For smoothly extended $\beta$ and $\phi$, 
we differentiate the boundary condition and obtain on 
$\partial\Omega$
\begin{equation}\label{int chi}
\frac{2\langle\nu,e_1\rangle}{\langle\nu,\beta\rangle}
u_{\beta\tau}=\frac{2\langle\nu,e_1\rangle}
{\langle\nu,\beta\rangle}
\left(\phi_j\tau^j-\tau^j\beta^i_ju_i\right)=:
\chi(x,Du).
\end{equation}
On the boundary, we get
$$u_{11}=u_{\tau\tau}+\chi+
\frac{\langle\nu,e_1\rangle^2}{\langle\beta,\nu\rangle^2}
u_{\beta\beta}.$$
Since $\chi(x_w,\cdot)=0$, the function
$$\tilde w(x,t):=e^{\alpha(\psi-u)+\gamma|Du|^2}(u_{11}-\chi)$$
satisfies $\tilde w(x_w,t_w)=w(x_w,\tau,t_w)$, moreover, 
for all $x\in\partial\Omega$ and $t\in[0,T]$
\begin{eqnarray}
\tilde w(x,t)&=&e^{\alpha(\psi-u)+\gamma|Du|^2}(x,t)\left\{u_{\tau\tau}(x,t)+
\frac{\langle\nu,e_1\rangle^2}{\langle\beta,\nu\rangle^2}
u_{\beta\beta}(x,t)\right\}\nonumber\\
&\le&\left\{1-\langle\nu,e_1\rangle^2
\left(1-\frac{\left\vert\beta^T\right\vert^2}{\langle\beta,\nu\rangle^2}
\right)-\frac{2\langle\nu,e_1\rangle\left\langle\beta^T,e_1\right\rangle}
{\langle\beta,\nu\rangle}\right\}\tilde w(x_w,t_w)\nonumber\\
&&{}+c\langle\nu,e_1\rangle^2e^{\alpha(\psi-u)+\gamma|Du|^2}(x,t)
\nonumber\\
&\le&\left\{1+c\langle\nu,e_1\rangle^2-\frac{2\langle\nu,e_1\rangle
\left\langle\beta^T,e_1\right\rangle}{\langle\beta,\nu\rangle}
+c\frac{\langle\nu,e_1\rangle^2}{\max\limits_{
\genfrac{}{}{0pt}{}{\xi\in T_x\partial\Omega}{|\xi|=1}
}u_{\xi\xi}(x)}\right\}\tilde w(x_w,t_w)\nonumber\\
&\le&\left\{1+c_1\langle\nu,e_1\rangle^2-\frac{2\langle\nu,e_1\rangle
\left\langle\beta^T,e_1\right\rangle}{\langle\beta,\nu\rangle}
\right\}\tilde w(x_w,t_w).\label{corr inequ}
\end{eqnarray}
Now, Lemma \ref{intro lem} ensures that we can choose $c_1$ in the
last inequality independent of $\alpha$ and $\gamma$. 

We may assume that $c_1$ in (\ref{corr inequ})
is chosen sufficiently large and $\beta$ is sufficiently
close to $\nu$ such that the expression in the last curly
brackets in (\ref{corr inequ}) is bounded below by $\frac{1}{2}$.

We define on $\ol\Omega\times[0,T]$
$$W:=\frac{e^{\alpha(\psi-u)+\gamma|Du|^2}(u_{11}-\chi)}
{1+c_1\langle\nu,e_1\rangle^2-\frac{2\langle\nu,e_1\rangle
\left\langle\beta^T,e_1\right\rangle}{\langle\beta,\nu\rangle}}.$$
Assume that $W$
attains its maximum at $(x_W,t_W)$ and $t_W>0$.

First, we address the case $x_W\in\partial\Omega$.
Observe that $W(x_W,t_W)\le\tilde w(x_w,t_w)=W(x_w,t_w)$.
At $(x_w,t_w)$, we get
$W_\beta\le0$, which implies that
\begin{equation}\label{W max RP}
u_{11\beta}+\alpha\delta_0u_{11}
\le c(1+(1+\gamma)u_{11}),
\end{equation}
using $\delta_0$ from \eqref{int delta0}.
At $x_w$, keeping the notation of Lemma \ref{mixed lemma},
we differentiate the boundary condition 
$u_\beta=\phi$ twice in direction $e_1$.
The a priori estimates obtained so far, and
the fact that $D^2u$, restricted to tangential directions,
is diagonal, yield
\begin{equation}\label{diff bcon 2}
u_{\beta11}\ge-c-2\beta_1^1u_{11}-2\beta^n_1u_{n1}.
\end{equation}
Then, combining
(\ref{W max RP}) and (\ref{diff bcon 2}) implies
$$c(1+(1+\gamma)u_{11})\ge\alpha\delta_0u_{11}-c.$$
 For $\alpha=\alpha(\gamma)$ sufficiently 
large, we get an upper bound on $u_{11}(x_w,t_w)$. This completes
the $C^2$-estimates, if $W$ attains its maximum on
$\partial\Omega$.

Now, we consider the case that $W$ attains its maximum
at $(x_W,t_W)$, $x_W\in\Omega$. We use 
$$\Gamma=-\log\left(1
+c_1\langle\nu,e_1\rangle^2-\frac{2\langle\nu,e_1\rangle
\left\langle\beta^T,e_1\right\rangle}{\langle\beta,\nu\rangle}\right)$$
in the following calculations. $\Gamma$
is well-defined as the argument of the logarithm is bounded
below by a positive constant. Moreover, the
$C^2\left(\ol\Omega\right)$-norm of $\Gamma$ is uniformly 
bounded independent
of $\alpha$ and $\gamma$.
We use that 
$$\log W=\alpha\cdot(\psi-u)+\gamma\cdot|Du|^2 
+\log (u_{11}-\chi)+\Gamma$$
attains its maximum at $x_W$. Of course, we may assume that
$1\le(u_{11}-\chi)(x_W,t_W)$. 
At $(x_W,t_W)$, we get 
\begin{align*}
0\le\frac{\dot W}{W}=&\alpha\left(\dot\psi-\dot u\right)
+2\gamma u^k\dot u_k+\frac{\dot u_{11}-\frac{d}{dt}\chi}
{u_{11}-\chi},\\
0=\frac{W_i}{W}=&\alpha(\psi-u)_i+2\gamma u^ku_{ki}
+\frac{u_{11i}-D_i\chi}{u_{11}-\chi}+\Gamma_i,\\
\intertext{and, in the matrix sense,}
0\ge&\frac{W_{ij}}{W}-\frac{W_iW_j}{W^2}\\
=&\alpha(\psi-u)_{ij}+2\gamma u^k_ju_{ki}+2\gamma u^k
u_{kij}\\
&{}+\frac{u_{11ij}-D_{ij}\chi}{u_{11}-\chi}
-\frac{(u_{11i}-D_i\chi)(u_{11j}-D_j\chi)}{(u_{11}-\chi)^2}
+\Gamma_{ij},\\
\end{align*}
where we have used that $\Gamma$ is time-independent.
We use $D_\cdot$ and $\frac{d}{dt}$ to indicate 
that the chain rule has not yet been applied. 
In the rest of the section, we drop the argument, if
we evaluate at $(x_W,t_W)$. We get
$$0\ge u^{ij}(\log W)_{ij}-\dot W.$$
Estimates for the time derivatives of $\psi$ and $u$,
the strict convexity of $\psi$ (\ref{lambda fix}),
the fact that $\Gamma\in C^2$ with uniform bounds,
and the differentiated flow equation
(\ref{flow eqn}) yield
\begin{align}\label{eins}
0&\ge2\gamma\Delta u+\frac{1}{u_{11}-\chi}
\left(u^{ir}u^{js}u_{ij1}u_{rs1}\right)
-u^{ij}\frac{(u_{11i}-D_i\chi)(u_{11j}-D_j\chi)}
{(u_{11}-\chi)^2}\nonumber\\
&{}\quad+\frac{1}{u_{11}-\chi}
\left(\hat f_{p_i}u_{i11}-c
-c\cdot\left|D^2u\right|^2\right)
+2\gamma u^k\hat f_{p_i}u_{ik}\\
&{}\quad+\frac{1}{u_{11}-\chi}
\left(\frac{d}{dt}\chi-u^{ij}D_{ij}\chi\right)
-c(\alpha+\gamma)+(\alpha\lambda-c)\tru.\nonumber
\end{align}
Direct calculations and (\ref{flow eqn}) imply
$$\frac{d}{dt}\chi-u^{ij}D_{ij}\chi\ge
-c\left(1+\left|D^2u\right|+\tru\right).$$
We use $\frac{W_i}{W}=0$ to get
$$2\gamma u^k\hat f_{p_i}u_{ik}+\frac{\hat f_{p_i}u_{11i}}
{u_{11}-\chi}\ge-c\cdot\alpha-\frac{c\cdot\left(1+
\left|D^2u\right|\right)}{u_{11}-\chi}-c.$$
Now, these estimates are applied to (\ref{eins}).
Let  $\theta\in(0,\frac12)$ be a small constant, to be fixed later.
First, we assume that, still at $(x_W,t_W)$,
$$(1-\theta)u_{\eta\eta}
\equiv(1-\theta)\max\limits_{|\xi|=1}u_{\xi\xi}
\le(u_{11}-\chi).$$
Here, a direction $\eta$, $|\eta|=1$, is chosen which 
corresponds to a maximal eigenvalue. 
Schwarz's inequality gives
$$u^{ij}(u_{11i}-D_i\chi)(u_{11j}-D_j\chi)
\le(1+\theta)u^{ij}u_{11i}u_{11j}
+\frac{c}{\theta}u^{ij}D_i\chi D_j\chi.$$
 From the definition of $\eta$, we get
$$
u^{ir}u^{js}u_{ij1}u_{rs1}
\ge\frac{\max\limits_{|\xi|=1}u^{js}u_{\xi j1}u_{\xi s1}}
{u_{\eta\eta}}
\ge\frac{1-\theta}{u_{11}-\chi}u^{ij}u_{11i}u_{11j}.
$$
Using $\theta\le\frac{1}{2}$ and $\frac{W_i}{W}=0$,
we get
\begin{align}
& u^{ir}u^{js}u_{ij1}u_{rs1}
-u^{ij}\frac{1}{u_{11}-\chi}(u_{11i}-D_i\chi)(u_{11j}-D_j\chi)
\nonumber\\
\ge & u^{ir}u^{js}u_{ij1}u_{rs1}
-\left(1+\theta\right)
\frac{1}{u_{11}-\chi}u^{ij}u_{11i}u_{11j}
-\frac{c}{\theta(1-\theta)}\frac{1}{u_{\eta\eta}}u^{ij}
D_i\chi D_j\chi\nonumber\\
\ge&-\theta\frac{2}{u_{11}-\chi}u^{ij}u_{11i}u_{11j}
-\frac{2}{\theta}\frac{c}{u_{\eta\eta}}\left(
1+\tru+\left|D^2u\right|\right)\nonumber\\
\ge&-c\theta(u_{11}-\chi)
\left(\tru+\alpha^2\,\tru+\gamma^2\left|D^2u\right|\right)
\nonumber\\
&-\frac{1}{\theta}\frac{c}{u_{\eta\eta}}\left(
1+\tru+\left|D^2u\right|\right).\nonumber
\end{align}  
Combining this inequality with (\ref{eins}) gives
\begin{align*}
0&\ge\left(2\gamma-c\,\theta\gamma^2-\frac{c}{\theta(\Delta u)^2}
-c\right)\Delta u-c(1+\alpha+\gamma)\\
&\quad+\left(\alpha\lambda-c-c\,\theta\alpha^2
-\frac{c}{\theta(\Delta u)^2}-\frac{c}{\Delta u}\right)\tru.
\end{align*}
 We fix $\gamma$, $\alpha=\alpha(\gamma)$ 
sufficiently large, and finally
$\theta=\theta(\gamma,\alpha)$ sufficiently small.
This implies an upper bound
on $u_{11}$. 
Note that as before, first we have fixed $\gamma$
and then $\alpha$.

Now, it remains to consider the case
$$(1-\theta)u_{\eta\eta}\ge(u_{11}-\chi).$$
We assume that we are in the non-trivial situation
where
\begin{equation}\label{1-theta}
\left(1-\tfrac{\theta}{2}\right)u_{\eta\eta}
\ge u_{11}.
\end{equation}
Set
$$\ol\Omega\times S^{n-1}\times[0,T]
\ni(x,\xi,t)\mapsto\tilde W(x,\xi,t)
=\frac{e^{\alpha(\psi-u)+\gamma|Du|^2}(u_{\xi\xi}-\chi)}
{1+c_1\langle\nu,e_1\rangle^2-\frac{2\langle\nu,e_1\rangle
\left\langle\beta^T,e_1\right\rangle}{\langle\beta,\nu\rangle}},$$
where $\chi$ is introduced in (\ref{int chi}).
Assume that $\tilde W$ attains 
its maximum at a positive time $t_{\tilde W}$ 
at $x_{\tilde W}\in\ol\Omega$ for a direction $\xi\in S^{n-1}$.
Assume further, that $x_{\tilde W}\in\partial\Omega$. 
If $x_{\tilde W}\in\Omega$, a modification of the proof 
for the case, when $W$ attains its maximum in
$\Omega\times(0,T]$, implies a bound for the 
second spatial derivatives of $u$. 
Using a decomposition of $\xi$ as
in (\ref{dec xi}), we obtain for $\beta$ 
sufficiently $C^0$-close to $\nu$
$$|\tau(\xi)|^2\le 1+c\left\Vert\beta^T\right\Vert_{C^0}.$$
As a direct consequence of this decomposition, 
we see that
$$u_{\xi\xi}\le u_{\tau(\xi)\tau(\xi)}+c.$$
We apply \eqref{corr inequ} and \eqref{1-theta}
\begin{align*}
\tilde W(x_{\tilde W},\xi,t_{\tilde W})
&\le |\tau(\xi)|^2W(x_w,t_w)+c\\
&\le \left(1+c\left\Vert\beta^T\right\Vert_{C^0}\right)W(x_W,t_W)+c\\
&\le \left(1-\tfrac{\theta}{2}\right)
  \left(1+c\left\Vert\beta^T\right\Vert_{C^0}\right)
  \tilde W(x_W,\eta,t_W)+c\\
&\le \left(1-\tfrac{\theta}{2}\right)
  \left(1+c\left\Vert\beta^T\right\Vert_{C^0}\right)
  \tilde W(x_{\tilde W},\xi,t_{\tilde W})+c,
\end{align*}
where $c=c(\alpha,\gamma)$. 
 For $\beta^T$ sufficiently small, we obtain a uniform bound on 
$\tilde W(x_{\tilde W},\xi,t_{\tilde W})$.

\section{Longtime existence and convergence}

So far, we have obtained uniform estimates on $\dot u$, $Du$, and
$D^2u$ as long as a smooth solution exists.
For $t=0$, we enclose our initial value $u_0$ by 
translating solutions. The maximum principle implies,
that $u$ will stay between the translating solutions.
We obtain that
\begin{equation}\label{encl}
-c+v^\infty\cdot t\le u\le+c+v^\infty\cdot t.
\end{equation}
We apply H\"older estimates for the
second derivatives due to Evans, Krylov, and Safonov, as
well as Schauder estimates, see \cite{lieberman}.
Since \eqref{flow eqn} has  no explicit $u$-dependence, we get
longtime existence
with uniform bounds on all higher derivatives of $u$.

To show convergence to a translating solution, we consider
$w:=u-u^\infty$. The following argument 
is similar to \cite{awcvp1994, os2001trans}.
Using the mean value theorem, we see that $w$ satisfies
a parabolic flow equation of the form
$$\left\{\begin{array}{r@{~=~}ll}
\dot w & a^{ij}w_{ij}+b^iw_i & \mbox{in~}\Omega,\\
w_\beta & 0 & \mbox{on~}\partial\Omega.
\end{array}\right.$$ 
Thus, the strong maximum principle implies, that $w$ is
constant or its oscillation is strictly decreasing in time. In the
first case, $u$ is already a translating solution. In the
second case, we wish to exclude that the oscillation is
strictly decreasing but does not tend to zero. If 
the oscillation of $w(\cdot,t)$ tends to 
$\epsilon>0$ as $t\to\infty$, 
we consider for $t_i\to\infty$
$$u^i(x,t):=u(x,t+t_i)-v^\infty\cdot t_i.$$
We have uniform estimates in any $C^k$-norm for 
the derivatives of $u^i$, and locally (in time)
uniform bounds for the $C^0$-norm, see (\ref{encl}).
Hence, a subsequence of 
the functions $u^i$ converges locally (in time) uniformly in any
$C^k$-norm to a
solution $u^*$ of our flow equation (\ref{flow eqn}) that
exists for all time. 
As the oscillation of $w=u-u^\infty$ is monotone in $t$,
we see that the oscillation of $u^*-u^\infty$
is equal to $\epsilon$, independent of $t$. This is a
contradiction to the strong maximum principle. If the
oscillation of $w$ tends to zero, we see that $u$ 
converges to a translating solution in $C^0$ as $t\to\infty$.
 Adding a constant to the translating solution
$u^\infty$, we may assume that $u\to u^\infty$ 
in the $C^0$-norm as $t\to\infty$. 
Interpolation inequalities of the form
$$\Vert Dw\Vert^2_{C^0(\ol\Omega)}\le c(\Omega)\cdot
\Vert w\Vert_{C^0(\ol\Omega)}\cdot
\left(\left\Vert D^2w\right\Vert_{C^0(\ol\Omega)}
+\Vert Dw\Vert_{C^0(\ol\Omega)}\right)$$
for $w=u-u^\infty$ and its derivatives imply smooth
convergence. The proof of Theorem \ref{main thm}
is complete.

\appendix
\section{Prescribing Gau\ss{} curvature for entire graphs}

Here, we present an application of our previous existence result on bounded
domains to construct 
unbounded hypersurfaces  with
prescribed Gau\ss{} curvature.
Assuming that the hypersurface is given as an entire graph,
the problem is to
find a solution of
\begin{gather}\label{gauss eqn}
\frac{\det D^2u}{\left(1+|Du|^2\right)^{\frac{n+2}2}}=g(x).
\end{gather}
Observe that this equation fits in the context of the present paper,
cf.\ (\ref{ell eqn})
\[
v^\infty=\log\det D^2u-\log f(x,Du),
\]
by defining
\[
f(x,p)=g(x)/h(p)
\quad\text{with}\quad
h(p)=h(|p|)=\left(1+|p|^2\right)^{-\frac{n+2}2}
\]
and
looking for translating solutions with speed $v^\infty=0$.

Following an argument of Altschuler, Wu \cite{awcvp1994},
we will construct entire rotationally symmetric translating solutions from  
solutions on growing disk-type domains.
Using the graph of the lower half sphere with suitably chosen radius,
it is a direct 
consequence of the strong maximum principle 
that there cannot exist an entire strictly convex
translating solution for any constant value of the Gau\ss{} curvature.
In the rotationally symmetric case,
we have the following result.
\begin{theorem}\label{riem thm}
Let $g(x)=g(|x|)$ be positive, smooth and integrable. 
There exists an  entire strictly convex solution $u$ to 
(\ref{gauss eqn}) if and only if 
\[
v:=\log\frac{\int_{\R^n} h(|p|)dp}{\int_{\R^n} g(|x|)dx}\ge0.
\] 
The solution constructed here has a uniformly bounded gradient if  
and only if $v>0$.
\end{theorem}
We remark that the rotationally symmetric setting does also allow
for a proof by reducing the problem to an ordinary differential equation.
\begin{proof}
We will use that for given constants $R>0$, $\rho>0$, there is a 
unique, strictly convex  solution 
$(v^\infty,u^\infty)$ 
to the following problem
\begin{align*}
\left\{\begin{array}{r@{~=~}ll}
\displaystyle 
e^{v^\infty}&
\det D^2u^\infty\;\frac{h(|Du^\infty|)}{g(|x|)}\quad&\text{in }B_R(0),\\
\displaystyle 
u_\nu^\infty&-\rho&\text{on }\partial B_R(0),\\
u^\infty(0)&0.&
\end{array}\right.
\tag{$*_{R,\rho}$}
\end{align*}
This is a direct application of Theorem \ref{ell thm}.
Furthermore, the uniqueness  of the solution 
implies its rotational symmetry.
The solution can also be found by imposing a second boundary value condition.

 From the strict convexity, we deduce that $Du^\infty$ is a diffeomorphism
from $B_R(0)$ onto $B_\rho(0)$.
Integrating $(*_{R,\rho})$, we obtain 
\[
e^{v^\infty}\int_{B_R}g(|x|)dx=\int_{B_R}\det D^2u^\infty\; h(|Du^\infty|)dx
=
\int_{B_\rho}h(|p|)dp,
\]
which uniquely determines the speed 
\[
v^\infty=v^\infty(R,\rho):=\log 
\frac{\int_{B_\rho}h(|p|)dp}{\int_{B_R}g(|x|)dx}
\]
as a function of the parameters $R,\rho$.
Note that $v^\infty(R,\rho)$ is  strictly decreasing in $R$ and
strictly increasing in $\rho$
with $v^\infty(R,\rho)\to-\infty$ for $\rho\to0$.

1. Nonexistence for $v<0$:
We argue similarly as in the afore mentioned case of constant Gau\ss{}
curvature. Now, we replace lower half spheres by suitably constructed
solutions on finite domains with arbitrarily large gradients at the boundary.
Since  $v<0$,  there exists a unique $\hat R$ such that
$\int_{B_{\hat R}}g(|x|)dx=\int_{\R^n}h(|p|)dp$.
Assuming that
we have  an entire solution $u$ of (\ref{gauss eqn}),
there is a $\rho>0$ such that
$|Du(x)|<\rho/2$ for all $|x|<\hat R$.
Now, take the unique $R<\hat R$ satisfying
\[
\int_{B_R}g(|x|)dx=\int_{B_{\rho}}h(|p|)dp
\]
and consider the solution $(v^\infty,u^\infty)$ of  $(*_{R,\rho})$.
By definition, we know that $v^\infty=v^\infty(R,\rho)=0$, 
hence, $u^\infty$ 
solves equation (\ref{gauss eqn}) in $B_R(0)$. Furthermore,
the Neumann boundary condition and our choice of $\rho$ yield
$|Du^\infty|>|Du|$ in a neighborhood of $\partial B_R(0)$.
Thus, there is a translate  $u^\infty+m$, $m\in\R$, 
which is strictly greater than $u$.
Now we shift back until the graphs touch first at a point
$x\in B_R(0)$.
By the strong maximum principle, this is impossible
as $u$ and $u^\infty$ solve the same elliptic equation
in $B_R(0)$.

2. Existence for $v\ge0$:
We construct our solution by choosing a sequence of 
increasing radii $R_k$ tending to $\infty$.
By the monotonicity properties of the function $v^\infty(R,\rho)$,
we can find for each $R>0$  a unique $\rho_R$ such that  
$v^\infty(R,\rho_R)=0$.
We remark that  $\rho_R$ is an increasing sequence.
Again,  Theorem \ref{ell thm} gives a unique smooth
rotationally symmetric solution $u_R$ 
to $(*_{R,\rho_R})$, which is
defined on $B_{R}(0)$ and satisfies $v^\infty=0$.
Note that for fixed $R$, the speed $v^\infty$ and the 
normal derivative at the boundary,
$-\rho_R$, are uniquely related. Hence, the solution $u_R$ must coincide
on smaller balls $B_{R'}(0)$, $R'<R$, with the previous solutions 
$u_{R'}$ to $(*_{R',\rho_{R'}})$.
Therefore, as $R$ tends to $\infty$,
the sequence $\{u_R\}$ will converge uniformly on compact sets
to a limit $u$, defined on all of $\R^n$. 
Clearly, $u$ is a rotationally symmetric solution to (\ref{gauss eqn}). 
Observe that the sequence $\rho_R$ will diverge in the case $v=0$, whereas
it stays bounded for  $v>0$.
This proves the boundedness of $|Du|$ in the latter case.
\end{proof}

Proceeding as in the existence part of the proof,
we get easily that
non-integrable $g(x)$ also allow for solutions provided that
$h(p)$ is non-integrable too.

This observation can be extended to
the function $h(p)$ arising in the 
equation of prescribed Gau{\ss} curvature in
Minkowski space
\begin{gather}\label{mink gauss eqn}
\frac{\det D^2u}{\left(1-|Du|^2\right)^{\frac{n+2}2}}=g(x).
\end{gather}
Hence, $h(p)=h(|p|)=\left(1-|p|^2\right)^{-\frac{n+2}2}$,
which  is not integrable on $B_1(0)$.
\begin{theorem}
For all positive and smooth functions $g(x)=g(|x|)$,
there exists an entire strictly convex solution $u$ to 
(\ref{mink gauss eqn}) satisfying $|Du|<1$.
Moreover, for a solution constructed here, $|Du|\le 1-\epsilon$,
$\epsilon>0$, if and only if $g$ is integrable.
\end{theorem}
\begin{proof}
We proceed similarly as in the proof of Theorem \ref{riem thm}
and use the notation introduced there.
Here, the speed function $v^\infty(R,\rho)$ is only defined for $\rho<1$.
But still, $v^\infty(R,\rho)$ is  strictly decreasing in $R$ and
strictly increasing in $\rho$.
In addition, $v^\infty(R,\rho)\to-\infty$ for $\rho\to0$
and $v^\infty(R,\rho)\to\infty$ for $\rho\to1$.
As in part 2 of the proof of Theorem \ref{riem thm},
we can find for any $R>0$ a unique $\rho_R\in(0,1)$ 
satisfying $v^\infty(R,\rho_R)=0$.
Since $\rho_R<1$, we can choose a smooth function $\tilde h(p)$, defined 
on $\R^n$, such that  $h(p)=\tilde h(p)$ for all $|p|\le\rho_R$.
Again, Theorem \ref{ell thm}
gives a unique smooth convex
rotationally symmetric solution $u_R$ 
to $(*_{R,\rho_R})$ with $h$ replaced by $\tilde h$, which is
defined on $B_{R}(0)$ and satisfies $v^\infty=0$.
We remark that the convexity of $u$ implies that $|Du|\le\rho_R$ 
on $B_{R}(0)$.
Thus, $u$ also solves the elliptic equation, if
we replace $\tilde h$ by the original $h$.
Again, for $R>R'$, 
$u_R$ will coincide with solutions $u_{R'}$ obtained 
on smaller balls $B_{R'}(0)$.
Hence, for $R$ tending to infinity,
$u_R$ converges to  an entire solution $u$ of (\ref{mink gauss eqn}).
In the present case, the sequence $\rho_R$ stays 
uniformly bounded away from 1 if $\int_{\R^n} g(|x|)dx<\infty$,
whereas  $\rho_R$ converges to 1  if $g$ is non-integrable.
\end{proof}

In the general case without rotational symmetry, a theorem
corresponding to Theorem \ref{main thm} in Minkowski space
can be obtained easily from the techniques of this paper,
provided that there holds a uniform a priori bound of the
form $|Du|<1-\epsilon$, $\epsilon>0$.

\section{Illustrations}

To illustrate the convergence of solutions,
we investigate numerically the  flow equation
$$\left\{\begin{array}{r@{~=~}ll}
\dot u & \log\det D^2u &\mbox{in~}\Omega\times[0,\infty),\\
u_\nu(\cdot,t)&(u_0)_\nu 
&\mbox{on~}\partial\Omega,\,t>0,\\
u (\cdot,0)& u_0 &\mbox{in~}\Omega
\end{array}\right.$$
on the ellipsoidal domain
$\Omega=\left\{(x,y)\in\R^2: 1.1\cdot\left(x^2+(2y)^2\right)<1\right\}$,
where
$u_0(x,y)=1.5x^2+y^2-0.1y^4$. 

The numerical integration
has been carried out on a $200\times100$ grid corresponding
to $[-1,1]\times[-0.5,0.5]\in\R^2$. Let $\Omega_0$ consist
of all grid points contained in $\Omega$, and
$\partial\Omega_0$ denotes those grid points
not contained in $\Omega_0$ such that 
one of the nearest neighbors belongs to $\Omega_0$.
For simplicity, we keep the same notation for 
the discretized quantities.

We use an explicit scheme for time integration.
The boundary condition is implemented as follows: 
For all $x_0\in\partial\Omega$ 
let $y_0:=x_0+\nu(x_0)\cdot \tau_0$, where
$\nu(x_0)$ is the normalized negative gradient of $x^2+(2y)^2$
and 
$\tau_0=\inf\{\tau:x_0+\nu(x_0)\cdot \tau\in \mbox{convex hull}(\Omega_0)\}$.
We set
$u(x_0):=u_0(x_0)-u_0(y_0)-u(y_0)$.
Here, $u(y_0)$ is obtained by linear interpolation. 

\begin{figure}[htb]
\typeout{includes bild0.ps,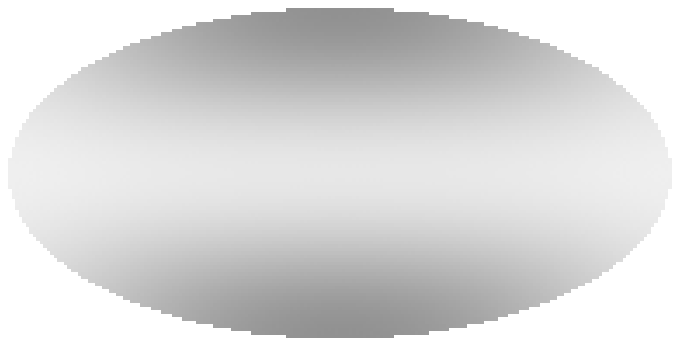,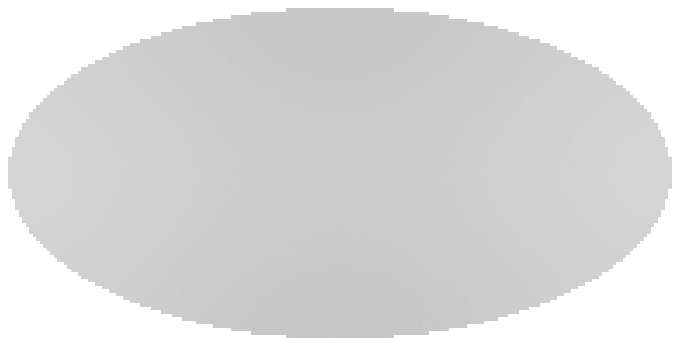}
\centering\mbox{%
\subfigure[$t=0.0$]
{\epsfxsize=0.3\textwidth\epsfbox{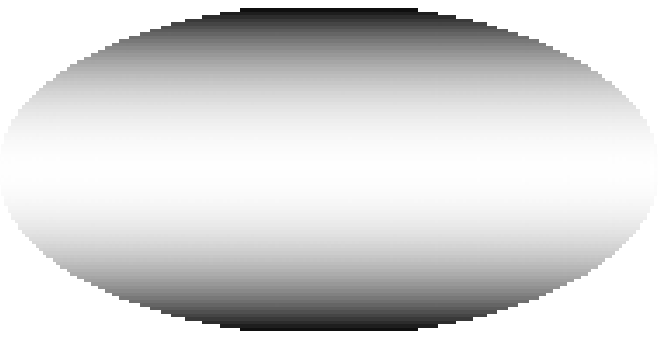}}
\subfigure[$t=0.1$]
{\epsfxsize=0.3\textwidth\epsfbox{bild1.ps}}
\subfigure[$t=0.5$]
{\epsfxsize=0.3\textwidth\epsfbox{bild5.ps}}}
\caption{Time evolution on an ellipsoidal domain} 
\label{bildchen}
\end{figure}
Figure \ref{bildchen}
shows a gray-scale plot of the velocity $\dot u$ at different times. 
It can be seen that the velocity tends to a constant,
reflecting the convergence of  $u$ to a translating solution.

\begin{figure}[htb]
\typeout{includes lpplot.ps + explot.ps}
\centering
\mbox{\subfigure[$\delta(t)$]
{\epsfxsize=0.5\textwidth\epsfbox{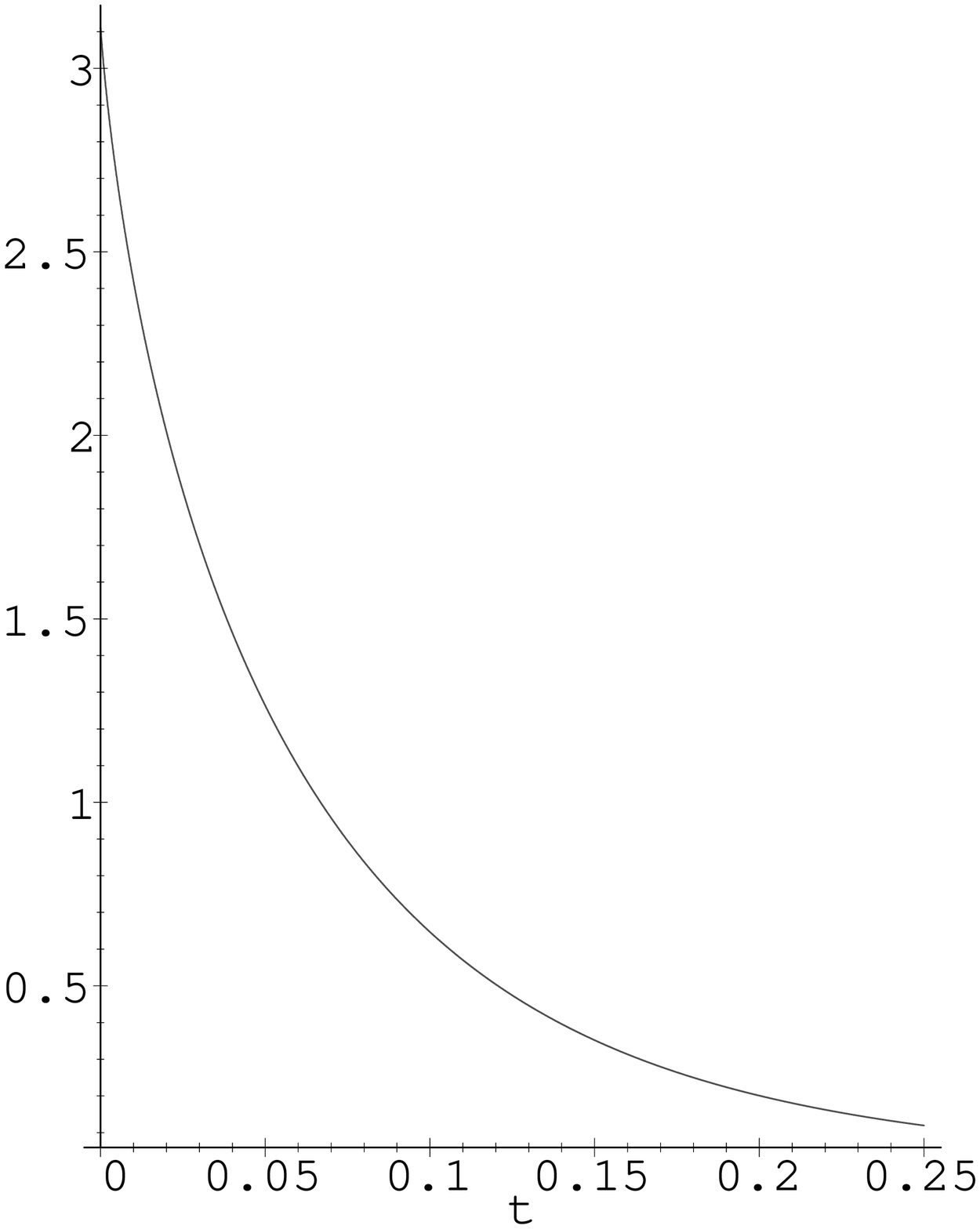}}
\subfigure[$\frac{-1}t\log \delta(t)$]
{\epsfxsize=0.5\textwidth\epsfbox{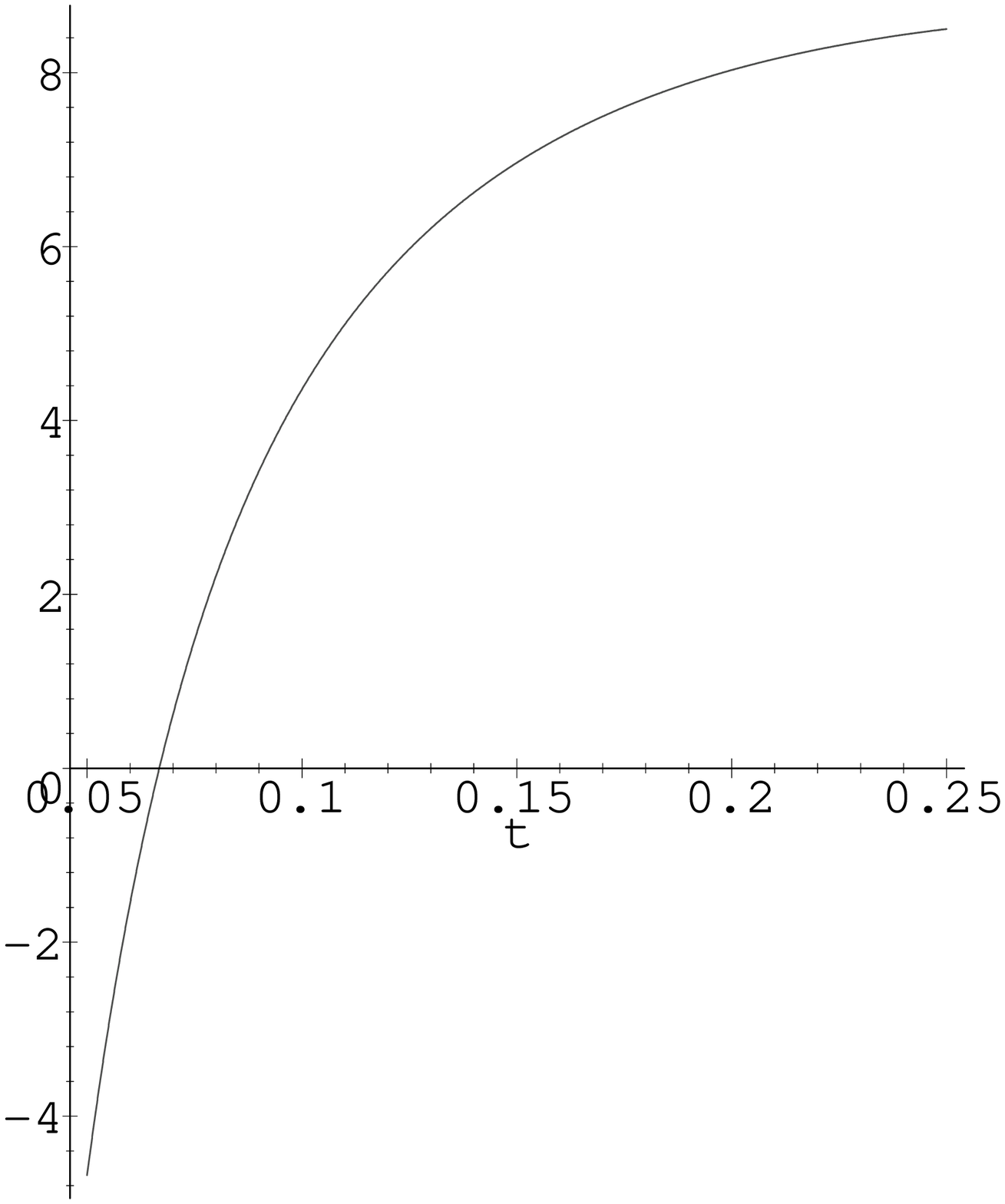}}}
\caption{Convergence to  constant velocity}
\label{bildchen2}
\end{figure}
In Figure \ref{bildchen2}a, we
show the decay of $\delta(t)=||\dot u(t)-\bar v(t)||^2_{L^2(\Omega)}$,
where $\bar v(t)=\frac1{|\Omega|}\int \dot u(x,t)\,dx$
is the mean velocity.
The expected exponential convergence can be seen from  
Figure \ref{bildchen2}b.
Here, we plot the exponential rate $\frac{-1}t\log \delta(t)$, 
which saturates 
nicely for larger times.

\bibliographystyle{amsplain}

\end{document}